\documentclass[10pt]{elsarticle}
\usepackage{amssymb}
\usepackage{amsmath}
\usepackage{amsfonts}
\usepackage{amsbsy}
\usepackage{amscd}
\usepackage{amsthm}
\usepackage{times}
\usepackage{psfrag}
\usepackage{graphics}
\usepackage{verbatim}
\usepackage{algorithm}
\usepackage{algorithmicx}
\usepackage{algpseudocode}
\usepackage{listings}
\usepackage{color}
\usepackage{colortbl}
\usepackage{graphicx}
\usepackage{caption}
\usepackage{subcaption}
\usepackage{array,supertabular}
\usepackage{tabularx}
\usepackage{booktabs}
\usepackage{multirow}
\usepackage{ulem}
\usepackage{units}
\usepackage{mathabx}
\usepackage{accents}
\usepackage{algorithmicx}
\usepackage{algorithm}
\usepackage{xfrac}
\usepackage{algpseudocode}
\usepackage{rotating}
\usepackage[capitalize]{cleveref}

\crefname{secinapp}{Section}{Sections}
\Crefname{secinapp}{Section}{Sections}

\newlength{\dhatheight}

\newtheorem{theorem}{Theorem}[section]

\newtheorem{definition}[theorem]{Definition}
\newtheorem{remark}[theorem]{Remark}

\def\XXint#1#2#3{{\setbox0=\hbox{$#1{#2#3}{\int}$ }
\vcenter{\hbox{$#2#3$ }}\kern-.57\wd0}}

\newcounter{myalgorithmctr}

\algnewcommand{\LineComment}[1]{\State \(\triangleright\) #1}

\newcommand{\func}[2]{\text{#1}\left(#2\right)}

\newcommand{\pluseq}{\mathrel{+}=}

\newcommand{\pr}[1]{\left(#1\right)}
\newcommand{\br}[1]{\left[#1\right]}
\newcommand{\set}[1]{\left\{#1\right\}}
\newcommand{\halfopen}[1]{\left[#1\right)}

\newcommand{\spanof}[1]{\text{span}\set{#1}}

\newcommand{\N}[2]{N_{#1}^{#2}}
\newcommand{\Nof}[3]{N_{#1}^{#2}\left(#3\right)}

\lstdefinestyle{FromFile}{language=Python,
                          frame=single,
                          numbers=left,
                          numberstyle=\tiny,
                          stepnumber=5,
                          numbersep=7pt,
                          numberfirstline=true,
                          abovecaptionskip=2pt,
                          belowcaptionskip=2pt}

\begin{document}

\begin{frontmatter}
\title{A Refinement Algorithm for Generalized B-splines}

\author[byu1]{Ian D. Henriksen}

\author[byu1]{Emily J. Evans\corref{cor1}}
\ead{ejevans@mathematics.byu.edu}

\cortext[cor1]{Corresponding author}

\address[byu1]{Department of Mathematics,
  Brigham Young University,
  Provo, Utah 84602, USA}

\begin{abstract}
In this paper we present a method for knot insertion and degree elevation of generalized B-splines (GB-splines) via the local representation of these curves as piecewise functions.  The use of local structures makes the refinement routines simpler to understand and implement.
\end{abstract}


\end{frontmatter}

\section{Introduction}\label{sec:intro}
In both design and analysis an important operation is that of refinement.  For a designer, local refinement allows curves to easily be modified so that the desired shapes can be created.  In analysis the techniques of both local and hierarchical refinement are used to increase numerical resolution in areas of the domain. Common types of refinement include $p$-refinement, $h$-refinement, $k$-refinement and $r$-refinement, however the most common refinement algorithms used in practice are degree elevation ($p$-refinement) and knot insertion ($h$-refinement). For B-splines, including non-uniform rational B-splines (NURBS), refinement algorithms are well studied and a wide variety of existing algorithms for refinement already exist~\cite{Fastdegelev,Prautzsch1984193,Boehm1980199,Piegl:1997:NB:265261}.   Refinement algorithms for generalized B-splines (GB-splines) are less well studied, in part, due to difficulty of the traditional definition of GB-splines using recursive integrals~\cite{trigspline, cheb}.  Although a method for knot insertion ($h$-refinement) is available for using recursive integrals~\cite{Ksasov1999, ue_spline_original, localgb}, the purpose of this paper to introduce new algorithms for GB-spline refinement using local structures.

Generalized GB-splines are a relatively new technology that seeks to overcome some of the shortcomings of NURBS. Rather than spanning the spaces of piecewise polynomials spanned by traditional B-spline curves, on each interval $\halfopen{t_i, t_{i+1}}$ in the given knot vector $T$, they span the spaces $\{1, t, \dots, t^{p-2}, u_i^{\br{p-1}}, v_i^{\br{p-1}}\}$ where $u_i^{\br{p-1}}$ and $v_i^{\br{p-1}}$ are $(p-1)^{th}$ integrals of arbitrary functions forming a Chebyshev space over $[t_i, t_{i+1}]$.
Because of their ability to span more general classes of functions, GB-splines allow exact representation of polynomial curves, helices and conic sections using control point representations that are intuitive and natural to designers~\cite{Ksasov1999,shapepreserv}.  Additionally GB-splines possess all of the fundamental properties of B-splines and NURBS that are important for design and analysis such as local linear independence, degree-elevation and partition of unity.  Another important aspect of GB-spines is that they behave similarly to B-splines with respect to differentiation and integration, thus facilitating the transfer of relevant properties of the continuous problem to the discrete problem~\cite{Ksasov1999,electromag}.

In 1999 Ksasov and Sattayatham \cite{Ksasov1999} demonstrated a variety of the properties of GB-splines.
In 2005 Costantini et al. \cite{costantini2005} studied generalized Bernstein bases of this form in greater detail.
In 2008, Wang et al. \cite{ue_spline_original} introduced unified extended splines (UE-splines), a subclass of GB-splines and demonstrated that this new class of splines contains several other classes of generalized splines.
In 2011, Manni et al. \cite{iga_gb_splines} proposed that GB-splines be used for isogeometric analysis.
In \cite{exponential_subdivision}, Romani successfully applied the techniques from \cite{subdivision_book} to form subdivision methods that allow for the approximation of UE-splines and other related classes of generalized splines via a limit of control meshes successively refined by a non-stationary subdivision scheme. In~\cite{quasiinterpolation}, quasi-interpolation was performed in isogeometric analysis using GB-splines.  GB-splines have also been used in the context of T-meshes~\cite{tmeshes}. Finally, in~\cite{gb-eval, insertinterestingnamehere} the authors introduced a new algorithm for the direct evaluation of GB-splines using local representations.  

In this paper an algorithm will be presented, using local representations of GB-splines that allows for knot insertion ($h$-refinement) and degree elevation ($p$-refinement) and follows a technique similar to the one used in \cite{projection_paper}.  The idea behind this algorithm is the creation of two knot vectors, $T_0$ which represents the original basis and $T_1$ which represents the refined basis.  The introduction of the local bases $V_i^p$ makes it so that any refinement operation of this form can be performed by computing the local representation for a spline over $T_0$ and then projecting it onto the desired basis functions over $T_1$.
Using local operations and projections to perform refinement adds additional efficiencies by allowing combinations of these knot insertions and degree elevations to be performed simultaneously.

\subsection{Structure and content of the paper}
\label{sec:content}
In~\cref{sec:GBRev} GB-splines are reviewed and appropriate notational
conventions are introduced.  In~\cref{sec:aegs} we also review an algorithm for the direct evaluation of GB-splines introduced earlier by the authors. Next in~\cref{sec:refine} we detail the refinement operations that can be performed by the routine. In~\cref{sec:proj} we present an algorithm for projection onto a refined basis including a description of the necessary variables and auxiliary routines. Next in~\cref{sec:aux} we present the auxiliary routines in greater detail.
Finally in~\cref{sec:conclusion} we draw conclusions.

\section{A review of generalized B-splines}
\label{sec:GBRev}
Generalized B-splines (GB-splines) were introduced in \cite{ue_spline_original}, and span spaces of the form $\set{1, t, \dots, t^{p-2}, u\pr{t}, v\pr{t}}$ where $u$ and $v$ are more general functions defined over each interval in a knot vector. The advantage of using GB-splines over traditional B-splines is that they allow for the exact representation of certain geometric curves and surfaces that cannot be well-represented by polynomial splines.  A more complete review of GB-splines along with their properties is given in the prior work by the authors~\cite{gb-eval} as well as the works by others including~\cite{Ksasov1999, ue_spline_original}.  We present a few necessary definitions for the algorithms presented later in this work.
\begin{definition} \label{knot_vector_definition}
A {knot vector}, $T$, is a nondecreasing vector of real numbers.  
\end{definition}
In this paper we will require all knot vectors to be open, in other words given a polynomial degree $p$ we require the first and last knots to have multiplicity $p+1$.
\begin{definition}\label{active region}
Given a spline of degree $p$, the active region of a given knot vector $T$ is the interval $[t_{p}, t_{m-p-1}]$. In short the active region of a given knot vector does not include the repeated end knots.
\end{definition}

\begin{definition}
Given a knot vector $T$, and functions $u_i$ and $v_i$ forming a Chebyshev space on each $\br{t_{i}, t_{i+1}}$ of nonzero length such that $u\pr{0} = v\pr{1} = 1$ and $v\pr{0} = u\pr{1} = 0$, we will refer to the sets of functions $u_i$ and $v_i$ as the {knot functions} over $T$.
\end{definition}
\begin{definition} \label{gbspline_definition}
Given a degree $p$ and a knot vector $T$ of length $m$ with corresponding knot functions $u_i$ and $v_i$.
Define the {$i^{th}$ GB-spline basis function of degree $p$}, denoted by $\N{i}{p}$ as follows:\\

\noindent Define the degree 1 GB-spline basis function as:
\[\Nof{i}{1}{t} = \begin{cases} u_{i}(t) & t \in \halfopen{t_i, t_{i+1}} \\ v_{i+1}(t) & t \in \br{t_{i+1}, t_{i+2}} \\ 0 & \text{otherwise} \end{cases},\]
For $p\geq 1$ define
\[\delta_i^p = \int_{t_i}^{t_{i+p+1}} \Nof{i}{p}{s} ds.\]
For $p\geq 1$ define $\Phi_i^p(t)$ as 
\[\Phi_i^p\pr{t} = \begin{cases} \frac{\int_{t_i}^t \Nof{i}{p}{s} ds}{\delta_i^p} & \text{if $\delta_i^p \neq 0$,}\\
0 & \text{if $ \delta_i^p = 0$ and $t < t_{i+p+1}$,}\\
1 &  \text{if $ \delta_i^p = 0$ and $t \geq t_{i+p+1}$.} \end{cases}\]
For $p > 1$, define
\[\Nof{i}{p}{t} = \Phi_i^{p-1}\pr{t} - \Phi_{i+1}^{p-1}\pr{t}.\]

In addition, if $t_{m-p-1} = \dots = t_{m-1}$ and $t_{m-p-2} \neq t_{m-p-1}$ (that is, if the last $p$ knots are repeated, and the last basis function is nonzero), define $N_{m-p-2}^p\pr{t_{m-1}} = 1$.
\end{definition}

\begin{definition}
Given a degree $p >1$, and a knot vector $T$ of length $m$ with a corresponding set of knot functions, a degree $p>1$, and $m-p-1$ control points $a_i$, define the corresponding {GB-spline curve} $f(t)$ as
\[f\pr{t} = \sum_{i=0}^{m-p-2} a_i \Nof{i}{p}{t}\]
for $t \in \br{t_{p}, t_{m-p-1}}$.
The curve $f(t)$ is not defined outside the active region of the knot vector $T$.
\end{definition}
\begin{remark}
In these definitions, and in the algorithms presented later in the manuscript we used zero-based indexing on the basis functions.  That is to say instead of indexing the basis functions from $1, \ldots, m$, we instead index the basis functions from $0, \ldots, m-1.$
\end{remark}

%
%
\begin{definition}
Given a knot vector $T$ with corresponding sets of knot functions $u_i$ and $v_i$, define
\[V_{i}^{p} = \spanof{1, \pr{t-t_i}, \dots, \pr{t-t_i}^{p-2}, u^{\br{p-1}}\pr{t}, v^{\br{p-1}}\pr{t}}\]
Where $u^{\br{p-1}}$, and $v^{\br{p-1}}$ are the $\pr{p-1}^{th}$ indefinite integrals of $u$ and $v$ respectively.
\end{definition}

\subsection{An algorithm for evaluation of GB-splines}
\label{sec:aegs}
Although Definition~\ref{gbspline_definition} is the traditional definition of a GB-spline, this definition does not provide for a simple means of evaluation. Using the definition the only effective means of evaluating GB-splines are either costly recursive numeric integration, or unwieldy symbolic computation of indefinite integrals. To address these difficulties a more direct method of computing values on GB-spline curves was presented in~\cite{gb-eval, insertinterestingnamehere}.  Since the ideas behind this algorithm are fundamental to the development of our algorithm for the refinement of GB-splines, especially the idea of local representations, we review the ideas of ~\cite{gb-eval} here.

The fundamental idea of the prior work was: since each basis function lies in the space $V_i^p$, a local representation of each basis function in terms of the functions spanning the space $V_i^p$ may be introduced.  This local representation was given by $u_i^{[p-1]}$, $v_i^{[p-1]}$, and an additional polynomial term of degree $p-2$ and was specifically chosen to be more amenable to integration.   With this local representation in mind, the definition of GB-splines in terms of their \textit{local representations} was given.
\begin{definition} \label{gbspline_standard_recurrence}
Given a degree $p$ and a knot vector $T$ of length $m$ with corresponding knot functions $u_i$ and $v_i$, and since $\N{i}{p}$ lies in the space $V_i^p$, we can represent $\N{i}{p}$ on the $j^{th}$ interval in $T$ as:
\[\Nof{i}{p}{t} = P_{i,j}^p(t) + a_{i,j}^p u_j^{[p-1]}(t) + b_{i,j}^p v_j^{[p-1]}(t)\]
where $P_{i,j}^p$ is a polynomial term and $a_{i,j}^p$ and $b_{i,j}^p$ are constants.
The recurrence stated in Definition \ref{gbspline_definition} can be written as:
\[a_{i,j}^p = \frac{a_{i,j}^{p-1}}{\delta_i^{p-1}} - \frac{a_{i+1,j}^{p-1}}{\delta_{i+1}^{p-1}},\]
\[b_{i,j}^p = \frac{b_{i,j}^{p-1}}{\delta_i^{p-1}} - \frac{b_{i+1,j}^{p-1}}{\delta_{i+1}^{p-1}},\]
and
\[\begin{aligned}P_{i,j}^p (t) &= \Nof{i}{p}{t_j} \\&+ \frac{1}{\delta_i^{p-1}} \left(\int_{t_j}^t P_{i,j}^{p-1}(s) ds - a_{i,j}^{p-1} u_j^{[p-1]}(t_j) - b_{i,j}^{p-1} v_j^{[p-1]}(t_j)\right) \\
	&- \frac{1}{\delta_{i+1}^{p-1}} \left(\int_{t_j}^t P_{i+1,j}^{p-1}(s) ds - a_{i+1,j}^{p-1} u_j^{[p-1]}(t_j) - b_{i+1,j}^{p-1} v_j^{[p-1]}(t_j) \right).
\end{aligned}\]

With the additional stipulation that if $\N{i}{p-1}$ is identically zero and the interval $[t_{i+p}, t_{i+p+1})$ is empty, then an additional 1 is added to $P_{i,j}$ to account for the modified treatment of basis functions that are identically $0$ in Definition \ref{gbspline_definition}.
As before, we also require that, if the last basis function is discontinuous at the end of the $m-p-1$ knot, that it must take a value of $1$ at that knot.
\end{definition}

\begin{figure}
\begin{minipage}[b]{0.5\linewidth}
\centering
\includegraphics[width=\textwidth]{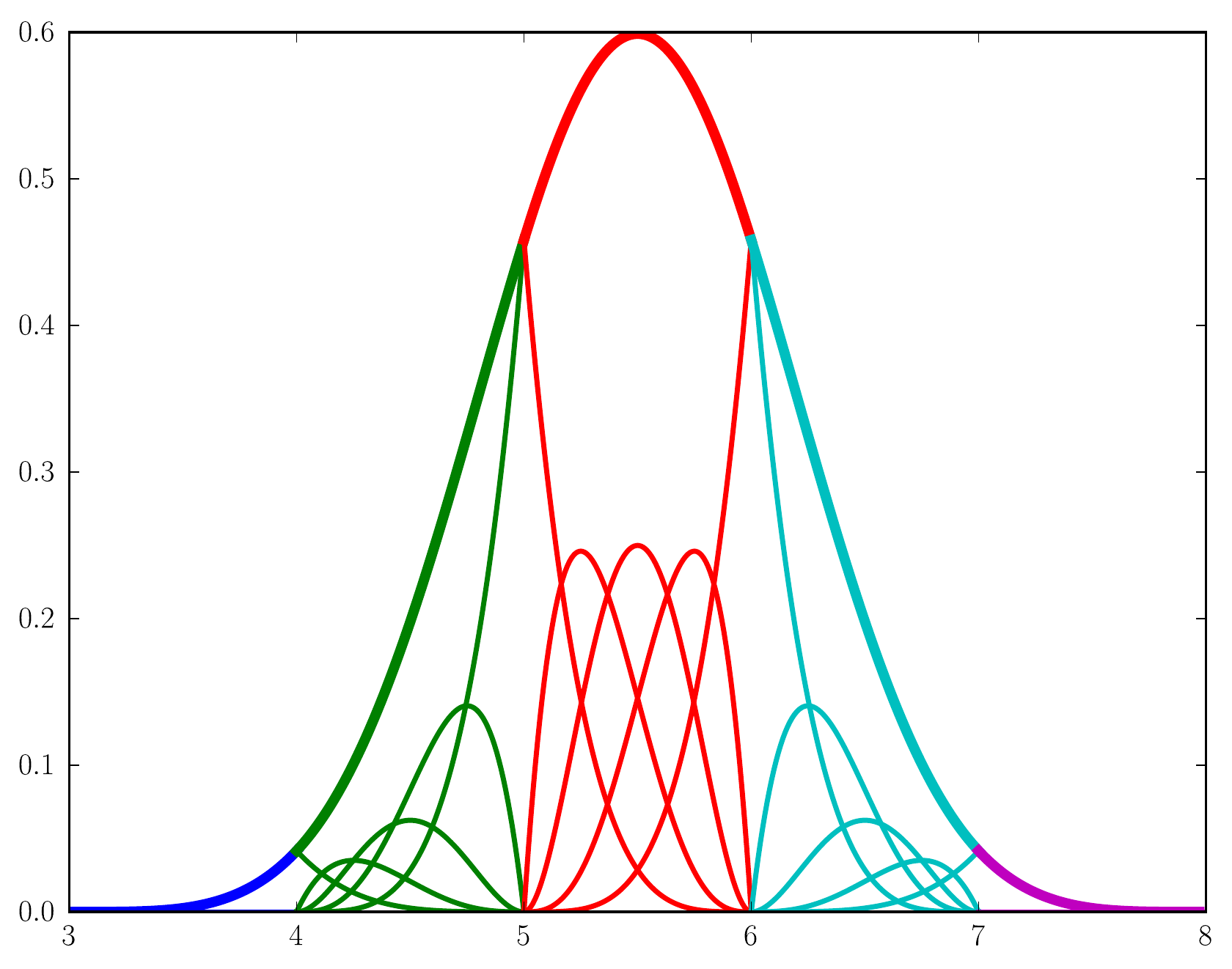}
\end{minipage}
\begin{minipage}[b]{0.5\linewidth}
\centering
\includegraphics[width=\textwidth]{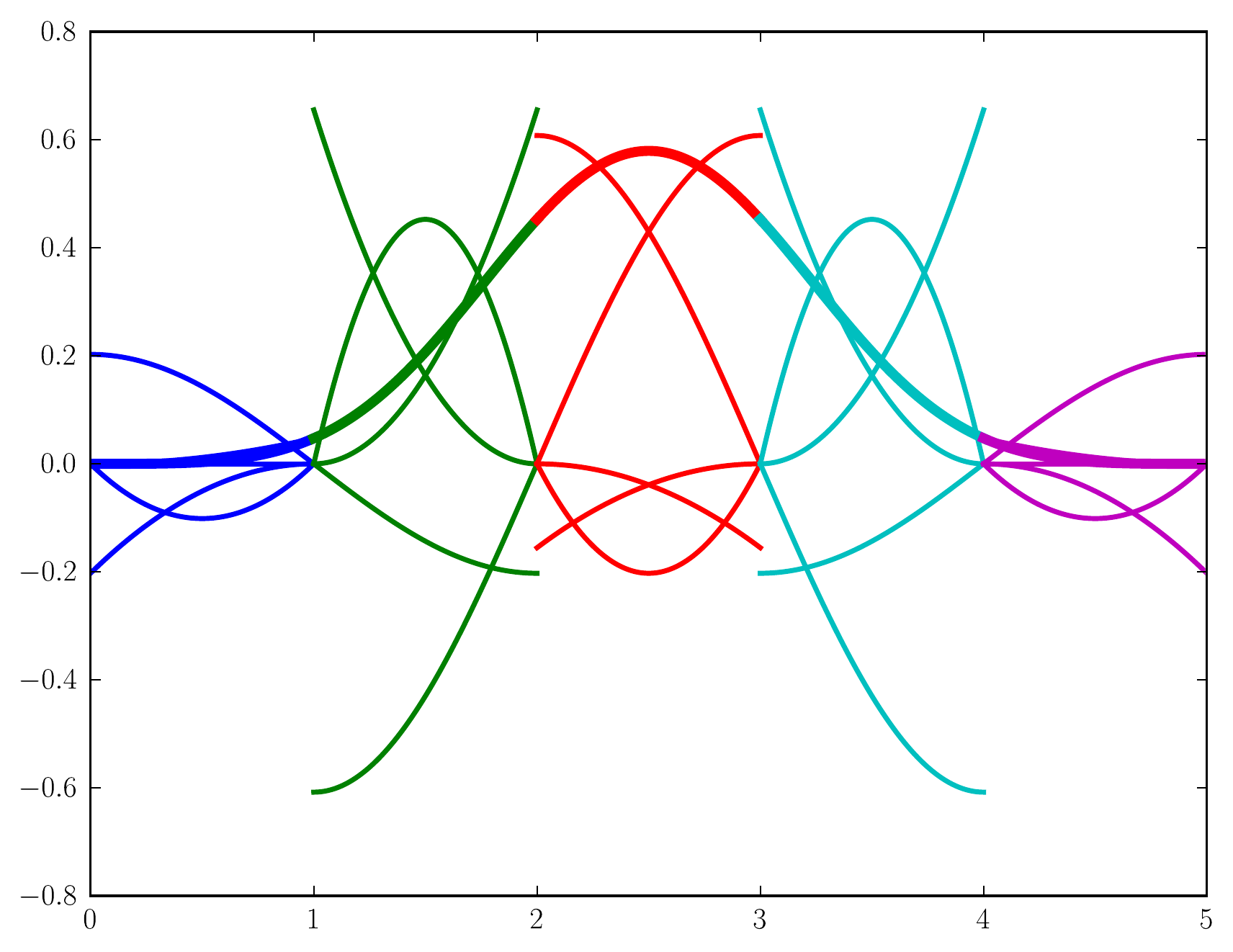}
\end{minipage}
\caption{The local polynomial representation for a uniform B-spline basis function of degree $4$ compared with the local representation for a uniform GB-spline function of degree $4$ defined using trigonometric knot functions.}
\end{figure}

The recurrence relation outlined in Definition \ref{gbspline_standard_recurrence} is not as easy to implement as De Boor's recurrence, however it does make it so that the evaluation of GB-spline curves is no longer tied to symbolic integrals or recursive quadrature.
It makes it clear that the values of $\N{i}{p}$ on the interval $[t_j, t_{j+1})$ depend only on $\Nof{i}{p}{t_j}$, the values of $u_j^{[p-1]}$ and $v_j^{[p-1]}$ at $t$ and $t_j$, and the full set of coefficients for $\N{i}{p-1}$ and $\N{i+1}{p-1}$.
These dependencies can be stated more explicitly.
To evaluate $\N{i}{p}$ at time $t \in [t_j, t_{j+1})$, it is necessary to know the values of the following function values:
\begin{itemize}
\item $u_j^{[p-1]}$ and $v_j^{[p-1]}$ at $t$;
\item the values of all the different $u_{j}, u_j^{[1]}, \dots, u_j^{[p-2]}$ and $v_j, v_j^{[1]}, \dots, v_j^{[p-2]}$ at $t_j$ and $t_j+1$ for each index $j$ corresponding to an interval in the support of $\N{i}{p}$;
\item the values of $u_j^{[p-1]}$ and $v_j^{[p-1]}$ if $j$ is an index corresponding to an interval in the support of $\N{i}{p}$ at $t_j$ if $t_j < t$  and at $t_{j+1}$ when $t_{j+1} < t$.
\end{itemize}

\section{Refinement operations on GB-splines}\label{sec:refine}
In design and analysis, given a spline curve the most common operation on the curve is refinement.  In design this refinement gives greater local control over the shape of the curve and in analysis this refinement gives greater numerical resolution.  Given a spline curve of degree $p$ represented by a knot vector $T_0$ we wish to represent the curve over a different knot vector $T_1$ where the active regions for $T_0$ and $T_1$ coincide and the spline basis $B_0$ over $T_0$ is contained in the span of the spline basis $B_1$ on $T_1$. 

If $T_1$ is the same as $T_0$ except that an additional knot has been inserted in the interior of the active region of $T_0$ then the operation under consideration is knot insertion.  If $T_1$ is formed from $T_0$ by inserting knots at each knot from $T_0$ that lies on the interior of the active region of $T_0$, adding an additional end condition knot on each side of the active region, and increasing the degree of the spline functions by 1, then the operation under consideration is degree elevation.  We note that the knots in the interior of the active region are repeated so that the continuity at the points in the knot vector is maintained in spite of the increase in the degree of the spline.

The introduction of the local bases $V_i^p$ makes it so that any refinement operation of this form can be performed by computing the local representation for a spline over $T_0$ and then projecting it onto the desired basis functions over $T_1$.
Refinement by projecting to and from local representations also provides a natural way to represent a spline over $T_0$ in terms of a knot vector $T_1$ with different end condition knots.
Projecting in this way also makes it so that any combination of these operations can be performed simultaneously.
As with the evaluation of GB-splines, the local piecewise representation of GB-splines can be used to avoid computing these recursive integrals.
In the case of refinement, projection to and from the local bases can also be used to perform degree elevation and change the placement of end condition knots.
The existence of a local basis also makes it so that coarsening operations like those performed on B-splines via B\'{e}zier projection in \cite{projection_paper} can be naturally extended to GB-splines, however in this paper we will restrict our attention to refinement only.

\begin{figure}[t]
\includegraphics[width=\textwidth]{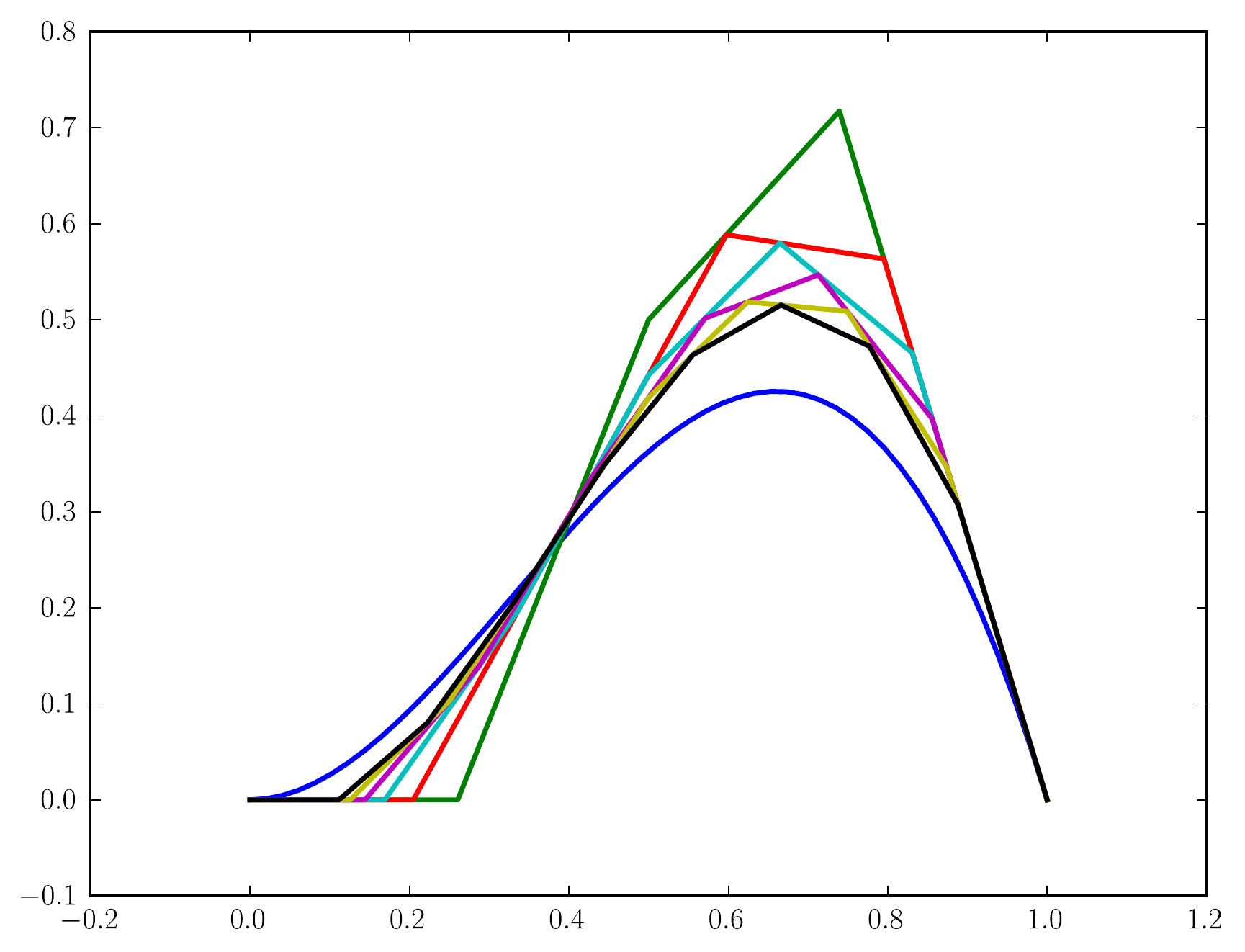}
\caption{Successive degree elevations of each basis function in a Bernstein-like basis of degree $3$ spanning trigonometric functions. The polygonal curves show the control meshes of successively higher degrees.}
\end{figure}

When refining GB-splines, we must also be certain that the integrals of the knot functions that are present in the space $V_i^p$ are spanned by the basis functions over the target knot vector.
For example, when performing degree elevation by a single degree, the knot functions $\tilde{u}$ and $\tilde{v}$ chosen over an interval $I_1$ in $T_1$ must be equal to a linear combination of the derivatives of the knot functions $u$ and $v$ over the interval $I_0$ in $T_0$ that contains $I_1$.
In general, it is sufficient for the knot functions $\tilde{u}$ and $\tilde{v}$ for an interval $I_1$ in $T_1$ used to form the target basis of degree $q$ over $T_1$ to be linear combinations of the $\pr{q-p}^{th}$ derivatives of the knot functions $u$ and $v$ over the interval $I_0$ in $T_0$ containing $I_1$.
These derivatives may not form a Chebyshev space over $I_1$, but this is sufficient.
For example, in the case of polynomial terms, any additional derivatives of the original $u$ and $v$ are no longer linearly independent of one another, and construction of the new basis fails.
On the other hand, as can be seen in \cite{projection_paper}, in the polynomial case, the polynomial terms over each basis provide a local basis that works perfectly well for refinement anyway.
Examples of when successive derivatives do yield a Chebyshev space include trigonometric and exponential splines.


When the derivatives of the knot functions cease to form a Chebyshev space, refinement is still possible as long as the span of $B_0$ is contained in the span of $B_1$.
All that is needed is a method for representing the integrals of the knot functions used in the local representations of $B_0$ in terms of the local bases used for $B_1$.

\section{An Algorithm For Projection onto a Refined Basis}\label{sec:proj}
Here we will present an algorithm that can be used to project a piecewise function in the span of a given spline basis $B_0$ defined on the active region of the knot vector $T_0$ corresponding to $B_0$ onto a different spline basis $B_1$, containing $B_0$ in its span, that is defined over a different knot vector $T_1$.
An additional routine will be provided to compute a piecewise representation from a set of control points and a set of local representations for a corresponding basis.


The main algorithm for projection onto a refined basis can be outlined as follows:
\begin{itemize}
\item Write the piecewise function in terms of the local basis used to construct the basis $B_1$.
	\begin{itemize}
	\item Subdivide the polynomial terms that correspond to intervals in $T_0$ that have been divided into new intervals in $T_1$.
	\item Degree-elevate the polynomial terms as many times as is necessary so that they are the same degree as the polynomial terms used for the target basis.
	\item Use the known values of the successive integrals of the knot functions to represent the knot functions over $T_0$ in terms of the local bases used for $T_1$.
	\item Add the polynomial terms from the representations of the knot functions over $T_0$ in terms of the local bases over $T_1$ into the values for the degree-elevated polynomials.
	\end{itemize}
\item Solve the linear systems corresponding to intervals with length that is neither zero nor almost zero as defined by a given tolerance to compute the coefficients for each basis function
	\begin{itemize}
	\item Aggregate the different computed coefficients for a given basis function by averaging them, raising an error if any value for an interval of non-negligible length is significantly different from the computed average.
	\end{itemize}
\end{itemize}
\begin{remark}
Aggregating the different computed coefficients by averaging them is not the only way to aggregate the coefficients, in fact using a weighted least squares method may give better results.  We use averaging for simplicity in implementation.
\end{remark}


This process is shown in greater detail in Algorithm \ref{gb_refinement}

\begin{remark}
The end condition knots are only needed for computing the local representations of either set of basis functions.
\end{remark}

\begin{algorithm}
\caption{Computing the control point representation for a given spline curve represented in piecewise form.}
\begin{algorithmic}[1]
\Procedure{RefineCurve}{\textit{polys0, genfunc0, treg0, rints0, polys1, genfunc1, reg1, rints1, tol}}
	\LineComment{Change the indexing of \textit{polys0} and \textit{genfunc0} to be indexed by the}
	\LineComment{intervals in \textit{regs1}.}
	\LineComment{Subdivide the polynomial terms so they are represented over the}
	\LineComment{intervals in \textit{regs1}.}
	\State \textit{polys0, genfunc0} $=$ RefineLocal(\textit{polys0, genfunc0, reg0, reg1, tol})
	\LineComment{Degree elevate the polynomials in \textit{polys0} so that}
	\LineComment{they are the same degree as the polynomials in \textit{polys1}.}
	\State \textit{polys0} $=$ ElevatePolys(\textit{polys0}, shape(\textit{poys1})[-1]-1)
	\LineComment{Compute the lengths of the intervals in $regs1$.}
	\State \textit{reg1lens} $=$ \textit{reg1}[1:] - \textit{reg1}[:-1]
	\LineComment{Find the intervals in \textit{reg1} that have positive length}
	\State \textit{pos} $=$ (\textit{reg1lens} $>$ \textit{tol})
	\LineComment{Represent the terms in \textit{genfunc0} in terms of the local basis on \textit{reg1}.}
	\State \textit{genfun0, offset} $=$ RepresentKnotFuncs(\textit{genfunc0, rints0, rints1, reg1lens, pos, tol})
	\LineComment{Add the polynomial part from the terms in \textit{genfunc0} to \textit{polys0}.}
	\LineComment{This makes it so that the function represented by \textit{polys0} and \textit{genfunc0} is now}
	\LineComment{completely rewritten in terms of the local bases used in \textit{polys1} and \textit{genfunc1}.}
	\State \textit{polys0} $\pluseq$ \textit{offset}
	\LineComment{Construct the local representations of the positive basis functions}
	\LineComment{for each interval in \textit{reg1}.}
	\State \textit{lbases} $=$ FullReverseDiagonals(concat(\textit{polys1, genfunc1}, $-1$))
	\LineComment{Combine the local representations in \textit{polys0} and \textit{genfunc0} into a single array.}
	\State \textit{lfunc} $=$ concat(\textit{polys0, genfunc0}, $-1$)
	\LineComment{Allocate \textit{coefs}.}
	\State $\text{\textit{coefs}} = $ empty array of \textit{NaN}'s of the same shape as \textit{lfunc}
	\LineComment{Find the coefficients that represent each piece of the piecewise function}
	\LineComment{as a linear combination of the basis functions in the target basis.}
	\State \textit{coefs}[\textit{pos}] = solve(swapaxes(\textit{lbases}[\textit{pos}, (0, 2, 1)), lfunc[\textit{pos},$\dots$,\textit{None}])$[\dots,0]$
	\LineComment{Combine the local results for each interval to compute the desired}
	\LineComment{basis coefficients.}
	\LineComment{Return the result.}
	\State \Return ReverseDiagonalAverages(\textit{coefs}, \textit{tol})
\EndProcedure
\end{algorithmic}
\label{gb_refinement}
\end{algorithm}

\subsection{Algorithm variables}
For ease in reading the algorithm we introduce the following variable naming convention and describe each variable.
Though they are not used directly in the algorithm, for the sake of explanation, we will let $p$ be the degree of the spline given in piecewise form and $q$ will be the degree of the target basis.
%
%
\begin{figure}[t]
\includegraphics[width=\textwidth]{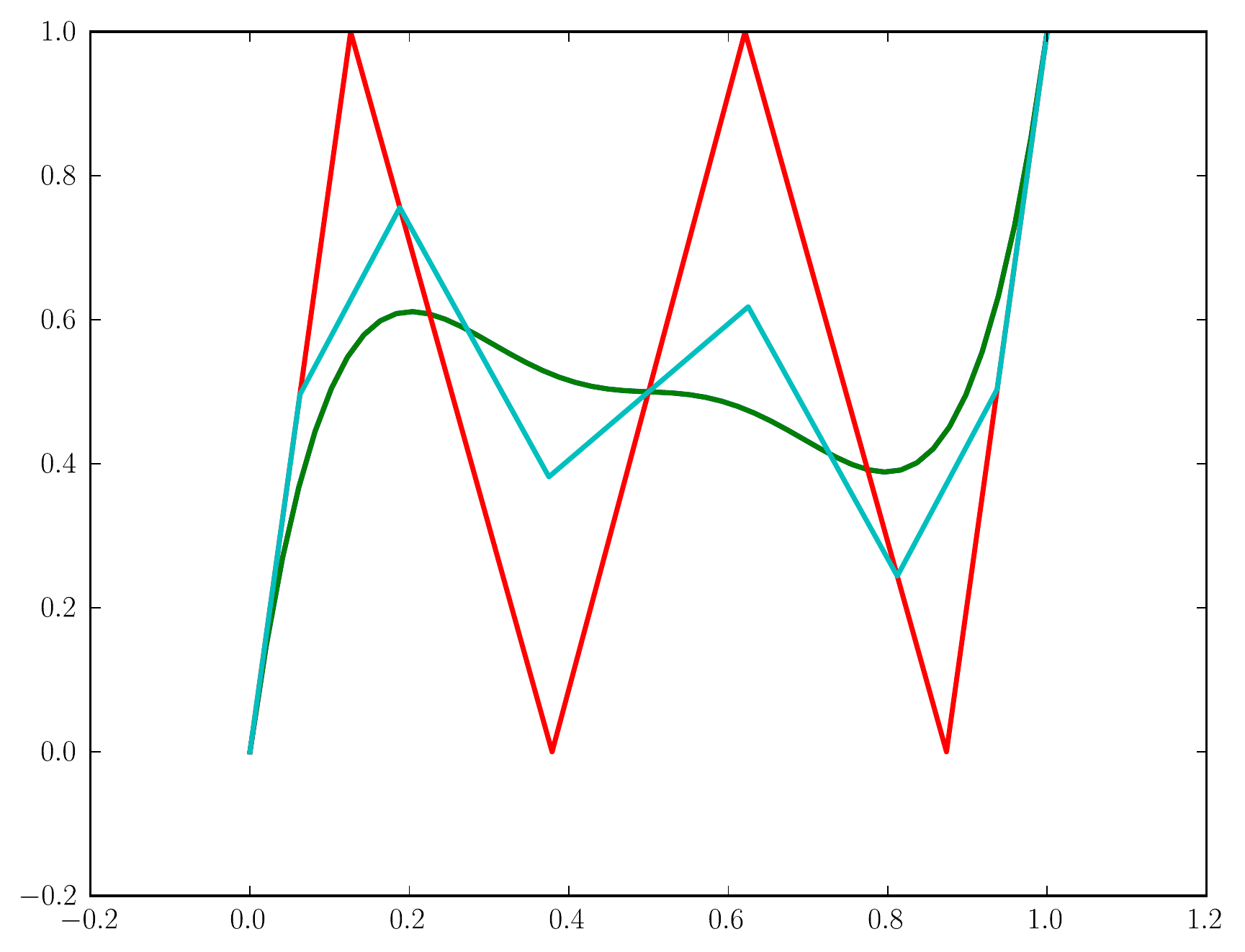}
\caption{Knot insertion at $.25$ and $.75$ on a degree $4$ GB-spline curve defined over $\br{0, 0, 0, 0, 0, .5, 1, 1, 1, 1, 1}$.}
\end{figure}

\begin{itemize}

\item \textit{tol} is a tolerance used to determine when intervals in a knot vector are short enough to be considered to have length $0$.
It will also be used to ensure the computed basis function coefficients do not differ significantly for the same basis function over different intervals.

\item \textit{polys0} are the polynomial terms for the piecewise function in the span of the first spline basis.
They are indexed first by interval within the active region of $T_0$, then by polynomial coefficient.
In this algorithm \textit{polys0} will be modified to contain the polynomial terms from the piecewise function represented over $T_1$.

\item \textit{polys1} are the polynomial terms for the basis functions in the target basis.
Basis functions are indexed along the first axis and intervals within the support of each basis function are indexed along the second axis.

\item \textit{genfunc0} is an array containing the coefficients for the general function terms of the knot functions corresponding to $T_0$.
It is indexed first by interval within the active region of $T_0$, then by the different general function terms.
The term for $u$ will be indexed at $0$ and the term for $v$ will be indexed at $1$.
 \textit{genfunc0} will be modified by the refinement algorithm so that it will represent the same general function terms except in the target basis.

\item \textit{genfunc1} is an array containing the coefficients for the general function terms of the knot functions over $T_1$ in the piecewise representation of the target basis.
Basis functions are indexed along the first axis, intervals within the support of each basis function along the second axis, and the two general function coefficients along the third (with $u$ first, then $v$).
This array has the shape of $\pr{n, p+1, 2}$.

\item \textit{reg0} is an array containing the knots that lie inside active region of $T_0$, with indexing and handling of knot repetition.
Regardless of whether or not the endpoints of the knot vector are repeated, if the piecewise term is a spline of degree $p$, this is an array equivalent to $T_0\br{p:-p}$.

\item \textit{reg1} is an array containing the knots that lie in the active region of $T_1$.  This array is similar to \textit{reg0} except that $T_1$ and the corresponding degree are used.

\item \textit{rints0} are the $0$ through $\pr{q-1}^{th}$ general function terms for the $\pr{q-p}^{th}$ derivatives of the knot functions over $T_0$ evaluated at the intervals in \textit{reg1}.
The number of integrals and the choice of integrals is dependent on $T_1$, but the knot functions are still chosen from $T_0$.
The $p^{th}$ integrals are indexed in ascending order along the first axis.
The different intervals within the knot vector are indexed along the second axis.
The different endpoints of each interval are indexed along the third axis.
The different knot functions ($u$ and $v$) are indexed along the last (fourth) axis with $u$ first.

\item \textit{rints1} is the $0$ through $\pr{q-1}^{th}$ general function terms for the knot functions over $T_1$ evaluated at the intervals in \textit{reg1}. Indexing is similar to that of \textit{rints0}

\item \textit{npolys} is a temporary array used to restructure \textit{genfunc1}.

\item \textit{ngenfunc} is a temporary array used to restructure \textit{genfunc0}.

\item \textit{reg1lens} is an array containing the widths of the intervals in \textit{reg1}.

\item \textit{pos} a boolean array where each entry is true if the corresponding entry in \textit{reg1lens} is above a given tolerance. If \textit{pos} is false at a given location this means that the interval in the know vector is essentially zero.

\item \textit{offset} is a temporary array of polynomial terms of the same shape as \textit{polys0} once it is represented over the intervals in \textit{regs1}.

\item \textit{offset1} and \textit{offset2} are arrays the same shape as \textit{offset} that are used to compute \textit{offset}.

\item \textit{ivals} and \textit{nints} are temporary arrays shaped like \textit{rints0}.

\item \textit{diffs} is a temporary array with one less entry along the first axis.

\item \textit{lbases} is an array containing the local representations of the spline basis functions in the target basis for each interval in the active region of $T_1$.
It is indexed first by interval in the active region, then by basis function index (indexing from $0$ to $q$ only along the basis functions in that are nonzero over that interval), then by term in the local representation.
The terms are ordered first by the polynomial terms with the last two entries corresponding to the coefficients for the general function terms from the knot functions.

\item \textit{lfunc} is an array containing the local representation of the piecewise function that we are representing as a spline with respect to the given basis.
It is indexed first by basis function, then by term within the local representation with the terms in the local representation ordered the same way they are for \textit{lbases}.

\item \textit{coefs} is an array containing the computed coefficients for the spline basis functions as they are computed over each interval.
Rows in \textit{coefs} corresponding to near $0$ length intervals in \textit{reg1} are set to \textit{NaN}.

\item $i$ and $j$ will be indices use to loop through \textit{reg0} and \textit{reg1} respectively.

\item \textit{jidxs} will be a list of indices for \textit{reg1} corresponding to the left endpoints of intervals of nonzero length that all lie within a given interval in \textit{reg0}.

\item \textit{endpts} an internal list of entries from \textit{reg1} corresponding to all endpoints of the intervals that have an interval indexed by an entry in \textit{jidxs}.
Shared endpoints between intervals are not repeated here.

\end{itemize}
\subsection{Auxiliary Routines}
In addition, when performing refinement, the auxiliary following routines are used:
\begin{itemize}

\item \textit{RefineLocal}: a routine that subdivides the polynomial terms in \textit{polys0} so that they are represented over the intervals of \textit{reg1}.
It also modifies \textit{genfunc1} so that it contains the coefficients for the integrals of the knot functions indexed by the intervals in \textit{reg1}.
The coefficients in \textit{genfunc1} are not refined by this routine, they are only copied into a new array that is indexed by the intervals of \textit{reg1} rather than by the intervals of \textit{reg0}.
The coefficients corresponding to an interval in \textit{reg1} are simply the coefficients corresponding to the interval in \textit{reg0} that contains it.

\item \textit{RestrictPoly}: a routine that restricts a polynomial over a given interval to a set of intervals and returns the new set of coefficients as an array.

\item \textit{ElevatePolys}: a routine to degree elevate an array of polynomial terms so they have a given degree.

\item \textit{LeftTaylorSeries}: a routine that computes the Taylor polynomial at the left endpoint of a given interval satisfying a given set of derivative values at that point.

\item \textit{RightTaylorSeries}: a routine that computes the Taylor polynomial at the right endpoint of a given interval satisfying a given set of derivative values at that point.

\item \textit{RepresentKnotFuncs}: a routine that changes the coefficients in \textit{genfunc0} from the representation for \textit{reg0} to the representation for \textit{reg1}.
This function also returns an array of polynomial terms to be added to \textit{polys1} to compensate for different allowable choices for integral values for the given functions.

\item \textit{FullReverseDiagonals}: a routine that takes a given array and, along the first two axes, forms a new array with rows that consist of the entries of all the reverse diagonals of full width.
Here it will always be called on arrays that have more rows than columns, so each row of the output array will have the same shape as a row of the original array.
This functions effectively as a method of translating between the piecewise representations of the basis functions and the variable \textit{lbases}.

\item \textit{ReverseDiagonalAverages}: a routine that takes the averages along all reverse diagonals of a given 2 dimensional array, ignoring $NaN$ entries.

\item
\textit{MatMul}: A utility function that performs matrix multiplication. Further details about this function can be found in~\cite{insertinterestingnamehere}.

\end{itemize}

\section{Auxiliary routines in detail}\label{sec:aux}
\begin{algorithm}
\caption{Restructuring \textit{polys0} and \textit{genfunc0} so they are indexed by interval in \textit{reg1}.
Subdivide the polynomial terms to represent them over the new intervals.}
\begin{algorithmic}[1]
\Procedure{RefineLocal}{\textit{polys0}, \textit{genfunc0}, \textit{reg0}, \textit{reg1}, \textit{tol}}
	\State \textit{npolys} = array of zeros of shape ( shape(\textit{reg1})[0] - 1, shape(\textit{polys})[1] )
	\State \textit{ngenfunc} = array of zeros of shape (shape(\textit{reg1})[0] - 1, 2 )
	\State $i = 0$
	\State $j = 0$
	\LineComment{Iterate through the intervals in \textit{reg0} and \textit{reg1}.}
	\While{$i < $shape(\textit{reg0})[0]}
		\LineComment{Skip the empty intervals in \textit{reg0}.}
		\While{\textit{reg0}[$i$] - \textit{reg1}[$i-1$] $< tol$}
			\State $i \pluseq 1$
		\EndWhile
		\LineComment{Stop if iteration has already reached the end of \textit{reg0}.}
		\If{$i = $ shape(\textit{reg0}[0])}
			\State \textbf{break}
		\EndIf
		\State \textit{jidxs} =  empty list
		\State \textit{endpts} =  empty list
		\LineComment{Iterate through the intervals in \textit{reg1} that lie in the current interval in \textit{reg0}.}
		\LineComment{Add the corresponding indices and endpoints to \textit{jidxs} and \textit{endpts}}
		\LineComment{if the interval has positive length.}
		\While{\textit{reg0}[$i$] - \textit{reg0}[$i-1$] $<$ \textit{tol}}
			\If{\textit{reg1}[$j+1$] - \textit{reg1}[$j$] $>$ \textit{tol}}
				\State append $j$ to \textit{jidxs}
				\State append \textit{reg1}[${j}$] to \textit{endpts}
			\EndIf
			\State $j \pluseq 1$
		\EndWhile
		\State append \textit{reg1}[${j}$] to \textit{endpts}
		\LineComment{Subdivide the polynomial term on the current interval onto the }
		\LineComment{corresponding intervals in \textit{reg1}.}
		\State \textit{npolys}[\textit{jidxs}] = RestrictPoly( \textit{polys0}[${i-1}$], \textit{reg0}[$i-1:i+1$], \textit{endpts})
		\LineComment{Copy the coefficients from \textit{genfunc0} into all corresponding }
		\LineComment{rows of \textit{ngenfunc}.}
		\State \textit{ngenfunc}[\textit{jidxs}] = \textit{genfunc0}[${i-1}$]
		\State $i \pluseq 1$
	\EndWhile
	\State \Return \textit{npolys}, \textit{ngenfunc}
\EndProcedure
\end{algorithmic}
\label{gb_refinement_RefineLocal}
\end{algorithm}
The auxiliary routine \textit{RefineLocal} begins the translation of the local representations over \textit{reg0} into local representations over \textit{reg1}.
It represents the polynomial terms over the new intervals.
It also restructures the coefficients in \textit{genfunc0} so that each row in \textit{genfunc0} corresponds to an interval in \textit{reg1}.
The coefficients stored in the row corresponding to a given interval $I_1$ are the coefficients that were originally stored for the interval $I_0$ in \textit{reg0} that contains $I_1$.
\textit{RefineLocal} performs these operations by iterating through the points in \textit{reg0} and \textit{reg1} and constructing the new arrays \textit{npolys} and \textit{ngenfunc} that contain the values that will be used to modify \textit{polys0} and \textit{genfunc0}.
This auxiliary routine is shown in greater detail in Algorithm \ref{gb_refinement_RefineLocal}.

The auxiliary routine \textit{RestrictPoly} is dependent on the polynomial representation used.
It accepts a polynomial term, the endpoints for the interval over which it is currently defined, and the endpoints for the intervals to subdivide the polynomial term onto.
It returns an array of polynomial terms containing the coefficients for the polynomial represented over each of the new intervals.
If the polynomial basis used is the Bernstein basis, this corresponds to subdividing a Bernstein polynomial.
If the polynomial basis used is the power basis shifted so that the starting point of each interval corresponds to $0$ in the polynomial term, this operation corresponds to a left shift operation. If the polynomial basis used is the power basis without any shifts, this operation does nothing.

The routine \textit{ElevatePolys} is also dependent on the polynomial representation used.
For the case of polynomials in Bernstein form, a simple algorithm follows easily from the fact that a polynomial in Bernstein form can be represented as a Bernstein polynomial of higher degree.
For power basis polynomials, this corresponds to filling in additional zero-valued coefficients to add higher degree terms in the polynomial representation of each set of coefficients.

The routine \textit{RepresentKnotFuncs} provides a way to transition between different choices that could possibly be made for the knot functions and their corresponding integral values.
First, since indefinite integrals are only well defined up to constant shifts, a method is needed to account for different choices for the general function terms.
Polynomial terms can be used to account for these differences in integral values.
In addition, if the $i^{th}$ set of knot functions are $\cos\pr{\frac{\pi}{2} t}$ and $\sin\pr{\frac{\pi}{2}t}$ over the interval $\br{0, 1}$ and the target intervals are $\br{0, \frac{1}{2}}$ and $\br{\frac{1}{2}, 1}$ with knot functions $\cos\pr{\frac{\pi}{2} t}$ and $\sin\pr{\frac{\pi}{2}t}$, and $\cos\pr{\frac{\pi}{2}\pr{t - \frac{1}{2}}}$ and $\sin\pr{\frac{\pi}{2}\pr{t - \frac{1}{2}}}$ respectively (taking linear combinations to satisfy the needed endpoint values), on $\br{\frac{1}{2}, 1}$, the knot functions over $\br{0, 1}$ are only equal to linear combinations of the knot functions on the interval $\br{\frac{1}{2}, 1}$.

On the other hand, since the polynomial terms used in the representation of a spline of degree $q$ have only degree $q-2$, any function in the span of $B_1$ is uniquely determined by its $0$ through $\pr{q-1}^{th}$ derivatives at each knot in the active region.
To see this, observe that, over a given interval $I_1$ with endpoints in \textit{reg1} and corresponding knot functions $\tilde{u}$ and $\tilde{v}$, any function $f$ on $I_1$ within the span of the basis functions can be represented uniquely in the form $a \tilde{u}^{\br{q-1}}\pr{t} + b \tilde{v}^{\br{q-1}}\pr{t} + P\pr{t}$.
Since the polynomial term $P\pr{t}$ has at most a degree of $q-1$, $f^{\pr{q-1}} = a \tilde{u} + b \tilde{v}$.
Since $\tilde{u}$ and $\tilde{v}$ form a Chebyshev space over $I_1$, they are linearly independent, so the coefficients $a$ and $b$ are uniquely determined.
Now observe that $P\pr{t} = f\pr{t} - a \tilde{u}\pr{t} - b \tilde{v}\pr{t}$, so the $0$ through $\pr{q-1}^{th}$ derivatives of $f$ at the endpoints of $I_1$ uniquely determine the $0$ through $\pr{q-1}^{th}$ derivatives of $P$.
Since $0$ through $\pr{p-2}^{th}$ derivatives at any given point uniquely determine a polynomial of degree $p-2$, $P$ is uniquely determined by its derivatives at either endpoint of $I_1$, so it is certainly uniquely determined by the $0$ through $\pr{q-1}^{th}$ derivatives at both endpoints.

The observation that the derivatives of a function at the endpoints of each interval uniquely define the local representation of the function also suggests the following method for computing the needed representations of the general function terms in the local representations of the original piecewise function.
Here, let $u_i$ and $v_i$ be the knot functions corresponding to the interval $I_0$ from \textit{reg0} that contains $I_1$.
Let $\tilde{u}_j$ and $\tilde{v}_j$ be the knot functions on $I_1$.
Also let $f_j$ be the restriction of the piecewise function being represented as a spline curve to the interval $I_1$.
It is already known that $f_j = g_j + P_j$ where $P_j$ is a polynomial term of degree $p-2$ and $g_j$ is a linear combination of $u_i^{\br{p-1}}$ and $v_i^{\br{p-1}}$.
\begin{itemize}

\item Compute the $0$ through $\pr{q-1}^{th}$ derivatives of $g_j$ at the endpoints of each interval with endpoints \textit{reg1}.

\item Represent $g_j^{\pr{q-1}}$ as $g_j^{\pr{q-1}} = a_j \tilde{u_j} + b_j \tilde{v_j}$.

\item Compute the $0$ through $\pr{q-2}^{th}$ derivatives of the general function terms $a_j \tilde{u}_j^{\br{q-1}} + b_j \tilde{v}_j^{\br{q-1}}$ at the endpoints of each interval in \textit{reg1}.

\item Use the computed derivative terms to find the $0$ through $\pr{q-2}^{th}$ derivatives of $g_j - a_j \tilde{u}_j^{\br{q-1}} - b_j \tilde{v}_j^{\br{q-1}}$.

\item Use the computed derivatives for the difference term to compute the polynomial part in the representation of $g_j$ in terms of the local bases used for $B_1$.

\end{itemize}
This process is shown in Algorithm \ref{gb_refinement_RepresentKnotFuncs}.

\begin{algorithm}
\caption{Restructuring \textit{polys0} and \textit{genfunc0} so they are indexed by interval in \textit{reg1}.
Subdivide the polynomial terms to represent them over the new intervals.}
\begin{algorithmic}[1]
\Procedure{RepresentKnotFuncs}{\textit{genfunc0, rints0, rints1, reg1lens, pos, tol}}
	\LineComment{Find the values of the $0$ through $(q-1)^{th}$ derivatives of each $g_j$.}
	\State \textit{ivals} = MatMul(\textit{rints0, genfunc0}$\br{\dots, None}$)
	\LineComment{Compute the representations of the $\pr{q-1}^{th}$ derivatives of $g_j$}
	\LineComment{in terms of the knot functions on the intervals in \textit{reg1}.}
	\State \textit{ngenfunc} =  an array of zeros the same shape as \textit{genfunc0}[$\dots$,\textit{New}]
	\State \textit{ngenfunc}[\textit{pos}] = solve(\textit{rints1}[0, \textit{pos}], \textit{ivals}[0, \textit{pos}])
	\LineComment{Compute the $0$ through $(p-1)^{th}$ derivatives of the general function terms}
	\LineComment{$a_j \tilde{u}_j^{\br{q-1}} + b_j \tilde{v}_j^{\br{q-1}}$.}
	\State \textit{nints} = MatMul(\textit{rints}[1\;:], \textit{ngenfunc})$\br{\dots, 0}$
	\LineComment{Take the difference between corresponding derivative terms.}
	\State \textit{diffs} = \textit{ivals}[1\;:, $\dots$, 0] - \textit{nints}
	\LineComment{Compute the polynomial term that satisfies the given derivative constraints.}
	\LineComment{Average the Taylor series at both endpoints to combine the results nicely.}
	\State \textit{offset0} = LeftTaylorSeries( reindex( \textit{diffs}$\br{::-1,:,0}, \pr{1, 0}$), \textit{reg1lens})
	\State \textit{offset1} = RightTaylorSeries( reindex( \textit{diffs}$\br{::-1,:,1}, \pr{1, 0}$), \textit{reg1lens})
	\State \textit{offset} $= .5 * \pr{\text{\textit{offset}}0 + \text{\textit{offset}}1}$
	\State \Return \textit{ngenfunc}$\br{\dots, 0}, \text{\textit{offset}}$
\EndProcedure
\end{algorithmic}
\label{gb_refinement_RepresentKnotFuncs}
\end{algorithm}

The routines \textit{LeftTaylorSeries} and \textit{RightTaylorSeries} both accept an array of derivative values, indexed first by interval, then by derivative degree from $0$ to $q-2$ and return a polynomial with the desired derivative values at the left and right endpoints of an interval of the given length with left endpoint at $0$.
\begin{remark}
The array of values passed to \textit{LeftTaylor Series} and \textit{RightTaylorSeries} creates and overdetermined system. If these two routines return values that are significantly different from one another, this usually indicates an ill-formed problem. If the output of these routines doesn't match, it may indicate, among other things, that the general functions for the piecewise function don't lie in the span of the general functions used for the target basis.
\end{remark}
The routine \textit{ReverseDiagonalAverages} is used to aggregate the results from the local computations over each interval in \textit{reg1}.
It takes a 2D array and computes the sum along each reverse diagonal ignoring entries of \textit{NaN}.
Since it should be aggregating terms that will be fairly close together, it can also be used to raise an error if any real entry differs from the average by more than a negligible amount. Errors here may also indicate ill-formed problems, specifically, bad values here may indicate that the piecewise function doesn't satisfy the continuity constraints built into the basis.
When it is called in Algorithm \ref{gb_refinement}, it is passed a tolerance for the purpose of identifying errors.

The helper routine \textit{FullReverseDiagonals} depends on the array library used.
Here it is only called on arrays with at least as many rows as columns, so it is sufficient to consider the cases where the length of a given reverse diagonal is equal to the length of a row of the original array.
%

To see how this indexing transformation works, let
\[a = \begin{bmatrix}
a_{0,0} & a_{0,1} & a_{0, 2} \\
a_{1,0} & a_{1,1} & a_{1, 2} \\
a_{2,0} & a_{2,1} & a_{2, 2} \\
a_{3,0} & a_{3,1} & a_{3, 2} \\
a_{4,0} & a_{4,1} & a_{4, 2} \\
a_{5,0} & a_{5,1} & a_{5, 2} \\
\end{bmatrix}\]
so
\[\func{FullReverseDiagonals}{A} = \begin{bmatrix}
a_{0,2} & a_{1,1} & a_{2,0} \\
a_{1,2} & a_{2,1} & a_{3,0} \\
a_{2,2} & a_{3,1} & a_{4,0} \\
a_{3,2} & a_{4,1} & a_{5,0} \\
\end{bmatrix}\]
For higher dimensional arrays, this indexing transformation is performed along the first two axes.
Conceptually, this can be thought of the same as the operation on $a$, except that the entries shown above are subarrays resulting from fixing two indices rather than single entries of a two dimensional array.
This indexing operation is used to transition between indexing local representations of basis functions first by basis function, then by interval within the support of each basis function, then by term in the representation, and indexing local representations first by interval in the active region of a spline, then by positive basis functions over that interval, then by term in the representation.

An additional benefit of the refinement algorithm presented here is that it provides a method for computing the Greville abscissae associated with a given spline basis.
The Greville abscissae are the coefficients for the linear combination of the spline basis functions that correspond to the linear function $y = x$.
These points are used to accurately represent $1$-dimensional spline curves as $2$-dimensional splines that are linear in one dimension.
This representation is commonly used to choose the $x$-axis values at which to plot the control points of a given $1$-dimensional spline curve.
In the case of $p=2$, these coefficients may not exist; however, if the spline basis spans degree $1$ polynomial terms, Algorithm \ref{gb_refinement} provides a method to compute them.

When performing degree elevation and knot insertion, it may also be necessary to rewrite a curve in piecewise form before passing it to the routine defined in Algorithm \ref{gb_refinement}.
This can be done very easily using the routine \textit{FullReverseDiagonals}, since it is essentially a transformation from indexing by basis function to indexing by interval.
The version presented here does involve some redundant computation, but, for the sake of simplicity, we will not optimize it further.
The process of forming a piecewise representation is shown in Algorithm \ref{FormPiecewise}.

\begin{algorithm}
\caption{Construct the piecewise representation of a spline curve given a set of control points and the corresponding set of basis functions.}
\begin{algorithmic}[1]
\Procedure{FormPiecewise}{\textit{cpts, polys, genfunc}}
	\LineComment{Construct the coefficients by multiplying the coefficients for each basis}
	\LineComment{function by the corresponding control points, reindexing so the coefficients}
	\LineComment{are accessed by interval, then summing the coefficients from the different}
	\LineComment{basis functions.}
	\State \textit{npolys} = sum( FullReverseDiagonals( \textit{polys} * \textit{cpts}[:,\textit{None,None}] ), 1)
	\State \textit{ngenfunc} = sum( FullReverseDiagonals( \textit{genfunc} * \textit{cpts}[:,\textit{None}, \textit{None}] ), 1)
	\State \Return \textit{npolys, ngenfunc}
\EndProcedure
\end{algorithmic}
\label{FormPiecewise}
\end{algorithm}
\section{Conclusion}
\label{sec:conclusion}
In this manuscript we have presented an algorithm for the refinement of GB-spline curves via their piecewise representation.  This algorithm was inspired by
previous work~\cite{gb-eval, insertinterestingnamehere} where an algorithm for the direct evaluation of GB-splines via their local representations was developed.  The use of 
local piecewise representations make refinement operations easier to understand and practical to perform.

The local structures used in this paper also provide insight into possible generalizations of the approach used in \cite{projection_paper} to local bases that are not generalized versions of the Bernstein basis.
Further generalizations of the approach used there would provide insights into possible structures that could be used for adaptive refinement of finite element meshes that satisfy given smoothness constraints on their boundaries.

In practice, GB-spline bases look increasingly similar to B-spline bases for successively higher degrees of spline bases.
A meaningful area for further work may come in explaining how, when, and how quickly this convergence occurs.
Providing meaningful bounds on this convergence would also make it much easier to determine when polynomial basis functions can be used as a replacement for GB-spline basis functions.

\bibliographystyle{elsarticle-num}
\bibliography{t2}

\end{document}